\documentclass[a4paper,11pt]{amsart}

\usepackage{amsmath}
\usepackage{amssymb,amsbsy,amsmath,amsfonts,amssymb,amscd}
\usepackage{latexsym}
\usepackage{geometry}
\usepackage{graphics}
\usepackage{color}
\input xy
\xyoption{all}

\input xy
\xyoption{all}

\newcommand\sE{{\mathcal E}}

\newcommand\sB{{\mathcal B}}

\newcommand\sR{{\mathcal R}}

\newcommand\Ga{\Gamma}
\newcommand\De{\Delta}

\newcommand\de{\delta}

\newcommand{\CC}{\ensuremath{\mathbb{C}}}

\newcommand{\ZZ}{\ensuremath{\mathbb{Z}}}

\newcommand{\hol}{\ensuremath{\mathcal{O}}}

\newcommand{\PP}{\ensuremath{\mathbb{P}}}

\newcommand{\ra}{\ensuremath{\rightarrow}}

\def\eea{\end{eqnarray*}}
\def\bea{\begin{eqnarray*}}

\newcommand\dual{\mathrel{\raise3pt\hbox{$\underline{\mathrm{\thinspace d
\thinspace}}$}}}
\newcommand\qe{\ifhmode\unskip\nobreak\fi\quad $\Box$}       

\def\BOX{\hfill\lower.5\baselineskip\hbox{$\Box$}}

\newtheorem{theo}{Theorem}[section]
\newtheorem{remarkk}[theo]{Remark}
\newenvironment{rem}{\begin{remarkk}\rm}{\end{remarkk}}

\newtheorem{defin}[theo]{Definition}
\newenvironment{definition}{\begin{defin}\rm}{\end{defin}}

\newtheorem{prop}[theo] {Proposition}
\newtheorem{cor}[theo]{Corollary}

\newtheorem{example}[theo]{Example}
\newenvironment{ex}{\begin{example}\rm}{\end{example}}

\newtheorem{question}[theo]{Question}
\newtheorem{fact}[theo]{Facts}

\newcommand{\sA}{\ensuremath{\mathcal{A}}}

\DeclareMathOperator{\Aut}{Aut}
\DeclareMathOperator{\Alb}{Alb}

\DeclareMathOperator{\Def}{Def}

\begin{document}

\title[Surfaces wit $p_g=0$]{Surfaces with $p_g=0$: constructions and 
moduli spaces, Burniat surfaces and deformations of
automorphisms}
\author{F. Catanese}

\thanks{The present article follows the expository thread of the talk 
given at the Kinosaki Algebraic Geometry Symposium,
October 27, 2011. It reports on research which  took place in the 
realm of the DFG
Forschergruppe 790 "Classification of algebraic surfaces and compact 
complex manifolds".}

\date{\today}
\maketitle
\section*{Introduction}

The following is a slightly  extended version of the talk, with the 
same title, which I gave at the Kinosaki Symposium on
Algebraic Geometry in October 2011, and
dealing with the classification of complex projective surfaces of 
general type (here the reader may find a few more references).

As mentioned in the talk, there exist two ways to do ``classification theory'':
\begin{itemize}
  \item one is similar to the activity of collecting beautiful and/or 
interesting objects at your home,
\item the other is like planning on the onset to build a large 
museum, starting by collecting financial support and experts
who are supposed to work there: in other words, giving  priority to 
the organizational and social aspects of your enterprise.
\end{itemize}

  But even if you choose method one, your home might become a museum 
after your death, so both methods
could converge to the same goal in the end; the main difference is 
therefore psychological, and the choice reflects
mainly  personal taste.

Also, at first glance the first method seems to be  simpler than the 
second: still it might
face you with non trivial technical problems (even if you only 
collect wine labels, it is not easy to peel them off the bottle,
especially for the French wines...).

Our beloved objects are here the surfaces of general type, and 
interesting patterns emerge while collecting
examples and studying them.

So, let $S$ be a smooth complex projective {\em minimal} surface of 
general type. This means that $S$ does not contain any
rational curve of self intersection
$(-1)$ or, equivalently, that the canonical divisor $K_S$ of
$S$ is nef and big ($K_S^2 > 0$). Then it is well known that
$$ K_S^2 \geq 1, \  \chi(S):= 1-q(S) + p_g(S) \geq 1.
$$

Recall that the {\em geometric genus} of $S$:
$$ p_g(S) :=h^0(S, \Omega^2_S) := \dim H^0(S, \Omega^2_S) = \dim 
H^0(S, \hol_S (K_S) ),
$$ and the {\em irregularity} of $S$:
$$ q(S) :=  h^1(S, \hol_S) := \dim H^1(S, \hol_S) = h^0(S, 
\Omega^1_S) := \dim H^0(S, \Omega^1_S),
$$ are {\em birational invariants} of $S$,
  as well as $K^2_S,$ since $$K^2_S = P_2 (S) - \chi (S), P_2 (S) : = 
\dim  H^0(S, \hol_S (2K_S) ).$$

Recall (see \cite{bombieri}) that the canonical model of $S$ is $X: = 
Proj (\sR (S))$, where $\sR (S): = \oplus_m H^0(S, \hol_S
(mK_S) )$ is the canonical ring of $S$.
$X$ is a normal surface with Rational Double Points as singularities, 
and with $K_X$ ample. Moreover there is a birational
morphism $ p \colon S \ra X$ contracting exactly the finitely many 
(-2)-curves (the irreducible curves $C$ such that
$ K_S \cdot C =0$, which satisfy  $ C^2 = -2$, hence  $C \cong \PP^1$ ).

  We have a coarse moduli space for the canonical models $X$ of 
surfaces $S$ of general type with fixed $\chi$ and $K^2$
(\cite{gieseker}).

\begin{theo} For each pair of natural numbers $(x,y)$  we have  the 
Gieseker moduli space $\mathfrak{M}_{(x,y)}^{can}$,
which is a quasi projective scheme, and
  whose points correspond to the isomorphism classes of minimal 
surfaces $S$ of general type with $\chi(S) = x$ and $K^2_S
=y$.

It is  a   coarse moduli space for the canonical models $X$ of 
minimal surfaces $S$ of general type with  $\chi(S) = x$ and
$K^2_S =y$.

\end{theo}

Concerning the range attained by the pair of numerical invariants 
above,  an upper  bound for $K_S^2$ is given by  the
Bogomolov-Miyaoka-Yau inequality:

\begin{theo}[\cite{miyaoka1}, \cite{yau},
\cite{yau2}, \cite{miyaoka2}] Let $S$ be a smooth surface of general type. Then
$$ K_S^2 \leq 9\chi(S),
$$ and equality holds if and only if the universal covering of $S$ is 
the complex ball $\mathbb{B}_2:=\{(z,w)
\in \CC^2 | |z|^2 + |w|^2 <1\}$.
\end{theo}

\subsection{Surfaces wih very low invariants}

  The above inequality is relevant when one is looking at the 
classification of surfaces of general type with `very
low' invariants, for instance with  the minimal possible value $\chi(S) =1$ for
$\chi (S)$.

In this case classification  means therefore to "understand" the nine 
moduli spaces $\mathfrak{M}_{(1,n)}^{can}$ for $1 \leq
n \leq 9$,
  in particular to describe their connected and irreducible components 
and their respective dimensions.
   Observe that
$$
\chi(S)=1 \ \ \iff p_g(S) = q(S).
$$
\begin{rem} If we assume that $S$ is irregular (i.e., $q(S) >0$), 
then by a result of Debarre (\cite{debarre}) it follows that
$K_S^2
\geq 2p_g(S)$.

Therefore in our case $p_g(S) \leq 4$. Moreover it was shown by 
Beauville (\cite{debarre})  that $p_g(S) = q(S) = 4$
if and only if $S$
is the product of two curves of genus 2.

  Also surfaces with $p_g=q=3$ were  described in \cite{ccml}: there 
are only two families, namely, the symmetric squares of
  curves of genus 3 (these have $K^2_S = 6$), or the quotients of a 
product  of curves $C_2 \times C'_3$,
  where $C_2 $ has genus $2$ and $C'_3$ has  genus $3$,  by
an involution of product type, acting as the hyperelliptic involution 
on $C_2 $  and freely on $C'_3$ (these surfaces have
$K^2_S = 8$).

In  \cite{ccml} a partial classification was shown (it was shown for 
instance that $6 \leq K^2_S \leq 9 $), and the
classification was then finished independently  by 
\cite{haconpardini}, and\cite{pirola}.

There has been lately a revival of interest for the surfaces in the 
second family, from the  topological side
(see \cite{akbulut}).

\end{rem}
The case $p_g= q=2$ is already harder: there is a substantial 
literature (see \cite{survey}), but they are still far from
being classified in spite of work by several authors: Zucconi, 
Ciliberto and Mendes Lopes, Hacon and Chen, Polizzi
and Penegini .

One knows that $4 \leq K^2_S \leq 9 $ by the cited inequalities, but 
it is unknown whether they do exist
for $ K^2_S =  9 $ (there have been repeated but failing attempts by 
Yeung to show that these ball quotients  cannot occur),
or whether the surfaces with $p_g= q=2$ and $ K^2_S =  4 $ are all 
double covers of a principally polarized Abelian
surface, and with branch curve a divisor $ \De \in | 2 \Theta |$ 
(this was proven by Manetti in \cite{manetti} under the
assumption that
$K_S$ is ample). There are also no examples known with $p_g= q=2, K_S^2 = 7$.

Surfaces with $p_g= q=1$ have $2 \leq K^2_S \leq 9 $ and have  been 
classified for $ K^2_S = 2,3 $ ($ K^2_S = 2
$ in \cite{cime0}, $ K^2_S = 3 $ in \cite{caci1}, \cite{caci2}, cf. 
also \cite{cp} for an exact  determination
  of the connected components
of the moduli space).  Existence of  the case of ball quotients 
($p_g= q=1, K^2_S = 9$) has been announced by Cartwright
and Steger.
  If this is correct,  there are no gaps and surfaces with  $p_g= q=1$ 
do exist   for each value of $K^2_S =2, \dots,   9$.

As the reader may know or surmise, the case
$p_g=q=0$ is the most difficult, and also the one for which there are 
more examples (see \cite{surveypg=0}). Surfaces
with $p_g=q=0$ exist for all values of $K^2_S =1, \dots,   9 $, but 
only the case $K^2_S =  9 $
has been classified by Cartwright and Steger, giving a precision to 
the fundamental work of Prasad-Yeung (\cite{ cs},
\cite{p-y}); there are 50 fundamental groups, and 100 isolated points 
of the moduli space, corresponding to 50 pairs of
complex conjugate non isomorphic surfaces (\cite{kk}).

Since the 1970's there was  a big revival of interest (see 
\cite{Dolgachev} for an early survey) in the
construction of these surfaces and in a possible attempt to
classification, and the Bloch conjecture and
differential topological questions raised by Donaldson Theory were a
further reason for raising further   interest
about surfaces of general type with
$p_g = 0$. There has been  recent  important progress
in the last 5 years (see  \cite{surveypg=0}) but there is no hope at 
the moment to even conjecturally finish the classification.
E.g., for $K_S^2 = 7$ there is only one family of surfaces with
$p_g=0$ known, constructed by Inoue
(cf. \cite{inoue}), while for $K_S^2 = 8$ the only known examples 
have the bidisk as universal cover (the reducible case has
been classified in \cite{ pg=q=0}, and a missing case was then added 
by \cite{ frap}).

At this point it seems appropriate to stop reporting on classification results
and to concentrate on `philosophical' issues. Consider  the following 
provocative question of D.
Mumford, posed at the Montreal Conference in 1980:

{\bf  Can a computer classify all surfaces of general type with $p_g = 0$?}

The meaning of the question is clear, and confirmed by recent 
progress: these surfaces are so many, that
it takes more than man's power to `conquer' their classification. And
it is indeed true that a computer algebra program is necessary to 
construct systematically certain surfaces, as it was carried
out in
\cite{4names}, \cite{bp} for product-quotient surfaces with $p_g=0$. 
More generally,
it is conceivable that computer programs, may be quantum computers, 
may describe all the possible canonical rings
of such surfaces in some not so distant future.

There remains however a major difficulty: these rings  will belong to 
different families, for instance according to the several
possibilities for the degrees of a minimal system of generators, and 
of relations
and higher syzygies.

But, how to find out how these locally closed subsets will fit in 
together inside the moduli space?

This difficulty is
witnessed already by the work of  Horikawa (\cite{horIII}) in the 
much simpler instance of surfaces with $p_g (S)= 4, K^2_S =
6$. Horikawa, looking at the canonical map, was able to divide these 
surfaces in 11
families, and began then to analyse the problem of incidence among 
these locally closed strata  of the moduli space,
the question being:
is stratum $\sA$  in the closure of stratum $\sB$?  This is a typical 
hard problem in the theory of surfaces,
and Horikawa showed that the corresponding subset of the moduli space 
has 4 irreducible components, and at most 3
connected components: leaving open the question whether there are 1, 
2, or 3 connected components.
Asnwering one of these questions turned out to be quite difficult, namely in
\cite{bcpAnnalen} it was shown that the number of connected 
components is at most 2, but it is still
open the question whether the number is 1 or 2.

One may sometimes be in a lucky situation, where it is possible  to 
describe completely a connected component of
the moduli space.

This may happen in several ways, for instance because there is only 
one mode of presentation for the canonical ring,
or because this phenomenon happens for some finite unramified 
covering $\hat{S}$ of $S$ ( see the next section).

Or, topology may dictate the existence of certain holomorphic maps to 
Abelian varieties or products of curves,
and this geometric feature  allows to  determine a connected 
component of the moduli space.

On the
other hand, if we are not in a lucky situation, or  if there is no 
good topological reason which determines a connected
component, it is very hard to show that an irreducible component is 
indeed a connected component. One has to study
deformations of a given family of surfaces (determining an open set 
$\mathfrak U$ in the moduli space), then one-parameter
limits of the deformed objects (degenerations of the surfaces 
corresponding to points in the open set $\mathfrak U$, i.e.,
determine the closure of the open set $\mathfrak U$) and then the 
deformations of these limits (try to see whether the
closure $\overline{\mathfrak U}$ is also an open set).

Together with Ingrid Bauer, also motivated by the problem of  sorting 
out the surfaces constructed in \cite{4names},
we took as a benchmark the problem of determining completely the 
connected components of the moduli spaces
containing the so called Burniat surfaces (some surfaces with $p_g 
(S) = 0$ constructed in 1966 by Pol Burniat, see
\cite{burniat}).  The problem is now solved for $K^2_S = 2,4,5,6$, 
but there is a single remaining
final step missing in the case
  $K^2_S = 3$ of tertiary Burniat surfaces.

The paper is organized as follows:

In the first section we shall briefly recall some by now classical 
`lucky' case where some connected component of
the moduli space of surfaces with
$p_g= 0$ can be determined. This part should also be seen as a `warm 
up' for the sequel.

In section 2 we define the Burniat surfaces and in section 3 we state 
the main classification theorem concerning them. In
section 4 we treat primary Burniat surfaces, which have a large 
fundamental group, and  we illustrate via  this case
the principle ``topology can determine connected components of the 
moduli space'', a phenomenon which has been
explored in various other cases.

In section 5 we introduce extended Burniat surfaces, which are 
deformations of nodal Burniat surfaces (they yield a concrete
example of an open set $\mathfrak U$ as previously mentioned).

Finally, in section 6 we describe a pathological behaviour of the 
moduli space, which is related to the degeneration of extended
Burniat surfaces to Burniat surfaces; namely, the fact that 
continuous families of canonical models yield,  at the level
of  minimal models,  families of branch loci which   vary discontinuously.
   The explanation goes through the remarkable phenomenon that, even 
if the automorphism group of the minimal model
is the same  as the automorphism group of the canonical model, the 
same does not hold for families;
so that, if $G$ is the group of automorphisms of the general surface, 
then  $\Def (S,G)$ is not proper onto $\Def (X,G)$;
and,  for tertiary Burniat
surfaces, while $\Def(S)$ and $ \Def(X,G)$ surject onto $\Def(X)$, $\Def 
(S,G)$ just maps to a nowhere dense set.

\begin{figure}[htbp]\label{configs}
\begin{center}
\scalebox{0.7}{\includegraphics{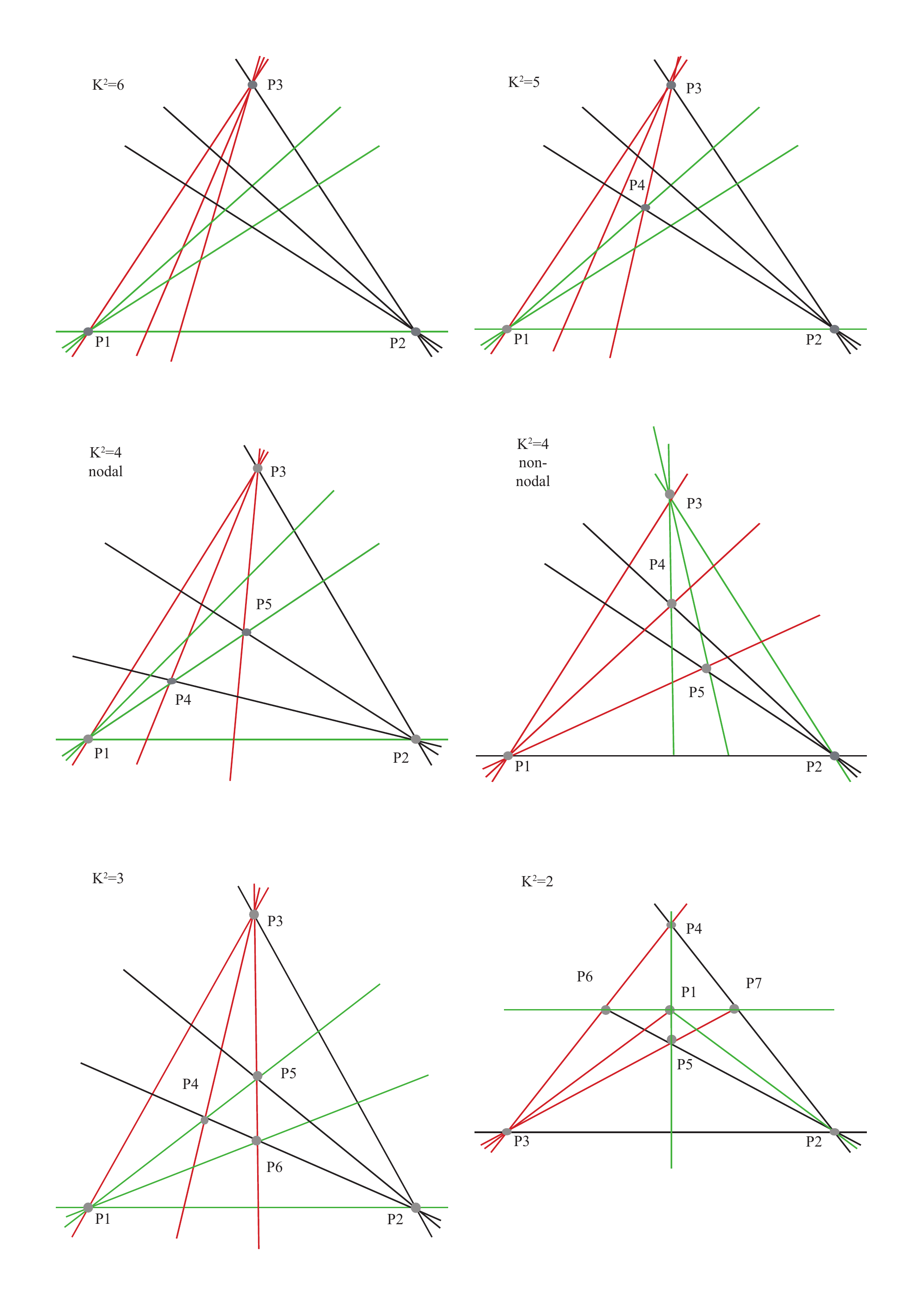}}
\end{center}
\caption{Configurations of lines}
\end{figure}

For the convenience of the reader we have drawn pictures of the line 
configurations in the plane corresponding to the branch
divisors of the Burniat surfaces with $K_S^2 = 2, \ldots 6$. They can 
be found in figure \ref{configs} attached below.

\section{Lucky cases}

Here are two classical examples of surfaces of general type where 
everything runs smoothly (see
\cite{miyaokagod}, \cite{miyaoka},\cite{tokyo},  \cite{Dolgachev}).

\subsection{Classical Godeaux surfaces}
\begin{enumerate}
\item
Here $ G : = \ZZ/ 5$, it acts on a 4-dimensional vector space $V$ via 
the 4 non trivial characters,
hence also on $\PP^3 : = \PP (V^{\vee})$ and $X =
\hat{X} / G$ is the quotient of an invariant 5-ic surface $\hat{X}$  (with RDP's 
as singularities) on which $G$ acts freely.

Hence, if $S$ is the minimal resolution of $X$,
$\pi_1 (S) \cong  G  = \ZZ/ 5$.
\item
It turns out that the representation $V$ is isomorphic to
the representation $H^0 (\hat{S}, \hol_{\hat{S}} (K_{\hat{S}}))$, 
therefore $S$ has $p_g(S) = 0$, and $K^2_S = 1$.
Since $G^{\vee} \cong   \ZZ/ 5$ is also the torsion part of $ H^2 (S, 
\ZZ)$, to $\chi \in G^{\vee} $ corresponds a divisor class
$M_{\chi}$, and
$$V \cong H^0 (\hat{S}, \hol_{\hat{S}} (K_{\hat{S}})) \cong [ 
\oplus_{\chi \in G^{\vee}, \chi \neq 0}H^0( S, \hol_S ( K_S +
M_{\chi}))]$$
\item
Conversely, if $S$ has $p_g(S) = 0$, and $K^2_S = 1$, and torsion $ T 
\cong G^{\vee} \cong \ZZ/ 5$, then the subspaces
$H^0( S, \hol_S ( K_S + M_{\chi}))$ have dimension 1, and they yield 
a basis for the vector space $H^0 (\hat{S},
\hol_{\hat{S}} (K_{\hat{S}}))$.

Let $ x_{\chi} \in H^0( S, \hol_S ( K_S + M_{\chi}))$ be a non zero 
element: then there cannot be any relation of the
form $ x_{\chi}  x_{\chi'}=  x_{\psi} x_{\psi'}$, because the 
associated divisors $div ( x_{\chi})$ are irreducible on the
canonical model $X$ of $S$. From this one concludes that the 
canonical map of $\hat{S}$ cannot have a quadric as image,
hence it induces an isomorphism of the canonical model $\hat{X}$ with 
a quintic surface in $\PP(V^{\vee}) = \PP^3$.
\end{enumerate}
\subsection{Standard Campedelli surfaces  }

\begin{enumerate}
\item
These are, by definition, the Campedelli surfaces with torsion group $T \cong ( \ZZ/ 2)^3$.
\item
Here   $\hat{S} $ is the natural unramified covering
associated to the torsion group (again here equal to the full first 
homology group $H_1 (S , \ZZ)$ since
$ q(S) = 0$), and  $S =
\hat{S} / G$.
\item
The best description of  $\hat{X} $ is as the maximal abelian 
covering of exponent 2 of the plane $\PP^2$
branched on 7 lines $ D_i$, one for each  $ g^{\vee}_i \in 
{G^{\vee}}^*:= G^{\vee} \setminus \{0\}$.
 $X$ is smooth if the 7 lines are in linear general
position.
\item
The Galois group of  $\hat{X} \ra \PP^2$ is the group $G'\cong (\ZZ/ 2)^6$
$$ G': = \oplus_{ g^{\vee}_i \in {G^{\vee}}^*}( \ZZ/ 2) g_i^{\vee} 
/(\ZZ/ 2) (\sum_i g_i^{\vee}) .$$
There is a natural surjection $ G' \ra G^{\vee}$, with kernel 
canonically isomorphic to $G$, since to each
element $ g \in G$ corresponds the sum of the elements lying in its 
annihilator $g^{\perp}= Ann (g)$ in $G^{\vee}$.

In this way we see that, being $X =
\hat{X} / G$,  $X\ra \PP^2$ is ramified on the 7 lines,
and  it has the property
that the local monodromy around the line
$D_i$ is the element $ g^{\vee}_i \in {G^{\vee}}$ in the Galois 
group. 
Instead $\hat{X} \ra X$ is unramified, with Galois group $G$.
\item

Indeed $\phi_{*} \hol_X= \hol_{\PP^2} \oplus \hol_{\PP^2} (-2)^7$, 
and $X$ is contained in the rank 7 vector bundle
$ \oplus_{ g_i \in G^*}L_i $ whose sheaf of sections is
isomorphic to $\hol_{\PP^2} (2)^7$. 

$X$ maps to  the fibre
product of the 7 double coverings
$$ y_{g_i}^2 = \Pi_ {g^{\vee}_j \notin  Ann (g_i)}  \delta_j ,$$
where $ D_j = div (\de_j)$,
and is indeed defined in the above rank 7 vector bundle by the following equations (see \cite{pardini}, and also \cite{cime},
page 146)
$$ y_{g_i} \cdot  y_{g_j}  =   y_{g_i + g_j} \Pi_ {g^{\vee}_j \notin ( Ann (g_i) \cup Ann (g_j))}  \delta_j ,$$
\item
Again we have a 7 dimensional vector space $V$ corresponding to  the 
7 non trivial characters of $G$,
$$V \cong H^0 (\hat{S}, \hol_{\hat{S}} (K_{\hat{S}})) \cong [ 
\oplus_{\chi \in G^{\vee}, \chi \neq 0}H^0( S, \hol_S ( K_S +
M_{\chi}))].$$
Each summand has dimension 1 and a generator $ x_{\chi}$ corresponds 
just to an equation $\delta_j$ for
a line $D_j$, using the established notation $\{  g^{\vee}_j\}$ for 
${G^{\vee}}^*$.

  The bicanonical map of $S$ is the Galois covering of $\PP^2$ with 
group $G^{\vee}$, and
$\hat{X} $ is embedded in $\PP^6 : = \PP (V^{\vee})$ as the complete 
intersection of 4 quadrics. The 4 quadrics are sums of squares, and are
easily obtained because
the 7 elements $ x_{\chi}^2$ belong to the 3 dimensional vector space 
$W \cong H^0( S, \hol_S ( 2 K_S ))=
p^* H^0( \PP^2 , \hol_{\PP^2 } ( 2  ))$.
\item
Conversely, given a Campedelli surface $S$   with torsion group $G 
\cong ( \ZZ/ 2)^3$, one considers the
natural unramified Galois covering   $\hat{S} \ra S$ with group $G$, 
and shows that $H^0( S, \hol_S ( K_S +
M_{\chi}))$ has dimension 1 for each character. Hence one has a 7 
dimensional vector space $V$
corresponding to  the 7 non trivial characters of $G$,
$$V \cong H^0 (\hat{S}, \hol_{\hat{S}} (K_{\hat{S}})) \cong [ 
\oplus_{\chi \in G^{\vee}, \chi \neq 0}H^0( S, \hol_S ( K_S +
M_{\chi}))].$$

Again, one can show that $\hat{X} $ is canonically embedded in $\PP^6 : = \PP (V^{\vee})$ as the complete 
intersection of 4 quadrics, and 
3 of the 7 elements $ x_{\chi}^2$
are linearly independent and yield the bicanonical map of $S$. Since the 4 quadrics are
sums of squares it follows that $\hat{X}$ is also invariant by the bigger group $G' \cong ( \ZZ/ 
2) ^6$, and then the bicanonical map of $S$
factors through the projection onto the canonical model $X$ and the 
Galois cover $X \ra \PP^2 $ with Galois group
$ G^{\vee} \cong G' / G.$

\end{enumerate}

\section{What is a ... Burniat surface?}
\noindent {\em Burniat surfaces} are surfaces of general type with 
geometric genus $p_g(S)=0$, and
they were constructed by Pol
Burniat in 1966 in
\cite{burniat}, where the method of singular bidouble covers was 
introduced in order to attack the geography problem for
surfaces of general type.

\noindent The birational structure of Burniat surfaces is rather 
simple to explain:

\noindent let $P_1, P_2, P_3 \in \mathbb{P}^2$ be three non collinear 
points (which we assume to be the points
$(1:0:0)$, $(0:1:0)$ and $(0:0:1)$), and let $D_i = \{ \Delta_i = 0 
\}$, for $ i \ \in \mathbb{Z} / 3 \mathbb{Z}$, be the union
of three distinct lines through $P_i$, including the line $D_{i,1}$ 
which is the side of the triangle joining the
  point $P_i$ with
$P_{i+1}$.

\noindent We furthermore assume that $ D = D_1 \cup D_2 \cup D_3$ 
consists of nine different lines.

\noindent
\begin{definition} A Burniat surface $S$ is  the minimal model for 
the function field $$ \mathbb{C} (\sqrt \frac{\Delta_1}
{\Delta_2}, \sqrt
\frac{\Delta_1} {\Delta_3} ).$$
\end{definition}

\begin{prop} Let $S$ be  a Burniat surface, and denote by $m$ the 
number of points, different from $P_1, P_2, P_3$, where
the curve $D$ has multiplicity at least 3 (hence indeed equal to 3). 
Then   $ 0 \leq m \leq 4$, and the invariants  of the
smooth projective surface
$S$ are:
$$ p_g(S) = q(S) = 0, \ K^2_S = 6 - m. $$
\end{prop} The heart  of the calculation, based on the theory of 
bidouble covers, as explained in \cite{sbc}, is that the
singularities where the three curves have multiplicities $(3,1,0)$ 
lower $K^2$ and the difference $p_g - q$ both by $1$, while
the singularities where the three curves have multiplicities 
$(1,1,1)$ lower $K^2$  by $1$ and leave $p_g - q$ unchanged
(in fact, for a bidouble cover branched on 3 smooth cubics, one has $ 
K^2_S = 9, p_g = 3$).

\begin{ex}(Singularities of Galois Coverings). Take three general 
lines $l_1,l_2,l_3$ through a point $P \in S$. Choosing
appropriate local (in the analytic topology) coordinates we can 
assume $P = 0 \in \CC^2$,  $l_1 = \{x=0\}$, $l_2 = \{y=0\}$
and $l_3 = \{x-y=0\}$. Taking the maximal Abelian cover of exponent 2 
branched in $l_1,l_2,l_3$ we get:
$$ u^2 = x, \ \ v^2 =y, \ \ w^2 = x-y \ \ \iff u^2 -v^2 =w^2,
$$

i.e. we get an ordinary double point
$$ Y:=\{(u,v,w) \in \CC^3 : u^2 -v^2=w^2\} \subset \CC^3,
$$ (an $A_1$-singularity) over $0$.

$Y$ is invariant under the involution
$$
\sigma \colon \CC^3 \rightarrow \CC^3, \ \ (u,v,w) \mapsto (-u,-v,-w),
$$

and we get a factorization of the $(\ZZ/2\ZZ)^3$- Galois covering 
$(Y,0) \rightarrow (\CC^2,0)$ as
$$ (Y,0) \rightarrow (Y/\sigma,0) \rightarrow (\CC^2,0).
$$ Note that $(Y/\sigma,0)$ is a $\frac 14(1,1)$ - singularity, which 
is not Gorenstein, but $\frac 12 $-Gorenstein.
Acquiring  such a singularity leaves the geometric genus $p_g$ and the irregularity $q$
invariant, but lowers $K_S^2$ by 1.

Indeed, the minimal resolution of such a singularity has an 
exceptional curve $E \cong \PP^1$ with $ E^2 = -4$, hence
the canonical divisor on the resolution $S$ is the pull back of the 
canonical divisor of $Y$ diminished by $\frac{1}{2} E$.
Because $ (K_S + E ) \cdot E = -2 \Rightarrow K_S  \cdot E = +2$.
\end{ex}

One may understand the biregular structure of a Burniat surface $S$ 
through the blow up $W$
of the plane at the points $P_1,
P_2, P_3,
\dots P_m$ of $D$  of multiplicity at least three.

$W$ is a  weak Del Pezzo surface of degree $6-m$ (i.e., a surface 
with nef and big  anticanonical divisor).

\begin{prop}\label{delpezzo} The   Burniat surface $S$  is a finite 
bidouble cover (a finite Galois cover with group
$(\mathbb{Z}/ 2
\mathbb{Z})^2$)  of  the weak Del Pezzo surface $W$.  Moreover   the 
bicanonical divisor $ 2 K_S$ is the pull back of the
anticanonical divisor
$ - K_W$. The bicanonical map of $S$  is the composition of the 
bidouble cover $ S \rightarrow W$ with the anticanonical
quasi-embedding of
$W$, as a surface of degree $ K^2_S = K^2_W$ in a projective space of 
dimension $ K^2_S = K^2_W$.
\end{prop}

\section{The main classification theorem} Fixing the number $K^2_S = 
6 - m$, one sees immediately that the Burniat
surfaces are parametrized by a rational family of dimension $K^2_S - 
2$, and that this family is irreducible except in the case
$K^2_S = 4$.

\noindent
\begin{definition} The family of Burniat surfaces with $K^2_S = 4$ of 
nodal type is the family where the points
$P_4, P_5$ are collinear with one of the other three points $P_1, 
P_2, P_3$, say  $P_1$.

\noindent The family of Burniat surfaces with $K^2_S = 4$ of 
non-nodal type is the family where the points
$P_4, P_5$ are never collinear with one of the other three points.
\end{definition}

Our main classification  result of  Burniat surfaces is summarized in 
the following table, giving information concerning the
families of Burniat surfaces, and where $\mathbb{H}_8 $ denotes the 
quaternion group of order 8. More information will be
given  in  the subsequent theorems.

  \begin{table}[ht]
  \small
\begin{tabular}{|c|c|c|c|c|}
\hline
$K^2$& dim  &  is conn. comp.? &name&$\pi_1$\\
\hline\hline

6&4&yes&primary& $1 \rightarrow \mathbb{Z}^6 \rightarrow \pi_1
\Rightarrow  (\mathbb{Z} /2 \mathbb{Z})^3$                   \\
\hline 5&3&yes&secondary& $\mathbb{H}_8 \oplus (\mathbb{Z} /2 \mathbb{Z})^3$ \\
\hline 4&2&yes&secondary& $\mathbb{H}_8 \oplus (\mathbb{Z} 
/2 \mathbb{Z})^2$ \\
 non nodal&&&&\\ 
\hline 4&2&no: $\subset$ 3-dim. irr.&secondary& $\mathbb{H}_8 
\oplus (\mathbb{Z} /2 \mathbb{Z})^2$ \\ nodal&&conn.
component $\supset$ &nodal& \\ &&$\supset$ extended Burniats&& \\
\hline 3&1&no: $\subset$ 4-dim. irr.&tertiary&$\mathbb{H}_8 \oplus 
\mathbb{Z} /2 \mathbb{Z}$\\ &&component $\supset$ &&  \\
&&$\supset$  extended Burniats&& \\
\hline 2&0& no: $\in$ conn. component&quaternary& $(\mathbb{Z} /2 
\mathbb{Z})^3$ \\ && of standard Campedelli&&\\
\hline
\hline
\end{tabular}
\end{table}

\noindent

\begin{theo} {\bf ( see  \cite{burniat1} and \cite{burniat2}) }

\ 1) The three respective subsets  of the  moduli spaces 
of minimal surfaces of general type $\mathfrak
M^{can}_{1, K^2}$   corresponding to Burniat surfaces with $K^2 = 6$, 
resp. with $K^2 = 5$, resp. Burniat surfaces with
$K^2 = 4$ of non nodal type, are irreducible connected components, 
normal, rational
   of respective dimensions 4,3,2.

\noindent Moreover, the base of the Kuranishi family of such surfaces 
$S$ is smooth.
\end{theo}

\noindent Observe that the above result for  $K^2 = 6$ was first 
proven by Mendes Lopes and Pardini in  \cite{mlp}. We
showed in
\cite{burniat1} the stronger theorem

\begin{theo}[Primary Burniat surfaces theorem]\label{primary} Any 
surface homotopy equivalent to a Burniat surface
with
$K^2 = 6$ is a Burniat surface with $K^2 = 6$.
\end{theo}

\begin{theo}[Secondary nodal Burniat surfaces 
theorem]\label{sec-nodal}  {\bf ( see  \cite{burniat2} and \cite{burniat3}) }

Secondary nodal Burniat
surfaces, together with secondary extended nodal Burniat surfaces 
form an irreducible connected component of the moduli
space.
\end{theo}

\noindent For $K^2 = 2$ another realization of the Burniat surface is 
(as shown by Kulikov in \cite{kulikov}) as  a special
element of the family of Campedelli surfaces with torsion 
$(\mathbb{Z} / 2 \mathbb{Z})^3$, considered in the previous
section.

  We saw that they are Galois covers  of the plane with group 
$(\mathbb{Z} / 2 \mathbb{Z})^3$, branched on seven lines.
For the Burniat surface we have the special configuration of a 
complete quadrilateral together with its three diagonals.

\section{The homotopy equivalence method} This section is devoted to 
the idea of the proof of theorem \ref{primary}. There
is a general philosophy behind the method of proof which applies to 
many more cases.

A Burniat surface $S$ with $K_S^2 = 6$ is called a {\em primary 
Burniat surface}. Recall that by proposition \ref{delpezzo}
$S$ is a finite bidouble cover of a Del Pezzo surface $Y = 
\hat{\PP}^2(P_1,P_2,P_3)$ of degree 6, which can be seen as
$$ Y:=\{((y_1,y'_1),(y_2,y'_2),(y_3,y'_3)) \in (\PP^1)^3 : y_1y_2y_3 
= y'_1y'_2y'_3\}.
$$

We take the $(\ZZ/2\ZZ)^3$-covering of
$$
\pi \colon \PP^1 \times \PP^1 \times \PP^1 \rightarrow \PP^1 \times 
\PP^1 \times \PP^1,
$$ given by
$$ (v_i)^2= y_i, \ \ (v'_i)^2 = y'_i, \ \ i \in \{1,2,3\}.
$$

Then $\pi^{-1}(Y)$ splits as the union of two Del Pezzo surfaces of 
degree 6, $Z \cup Z'$, where $Z:=\{v_1v_2v_3 =
v'_1v'_2v'_3\}$, and
$Z':=\{v_1v_2v_3 = - v'_1v'_2v'_3\}$. What we have done is the 
following:  we have taken the square root of the two
points in each
$\PP^1$ corresponding to two of the four lines passing through each 
$P_j$. Now we take the square root of the other two lines
through each of the three points $P_j$ and obtain a 
$(\ZZ/2\ZZ)^3$-covering $\sE_1 \times \sE_2 \times
\sE_3 \rightarrow \PP^1 \times \PP^1 \times \PP^1$, where each 
$\sE_j$ is therefore an
elliptic curve.  We get the following diagram:
\begin{equation*}
\xymatrix{
\hat{X} \ar[r]^{G} \ar[d]_{\hat{i}}& X' = \hat{X} /G\ar[dr] & \\
\hat{X } \cup \hat{X}' \ar[r] \ar[d] & Z \cup Z' \ar[d] \ar[r] &Y \ar[d]\\
\sE_1 \times \sE_2 \times \sE_3 \ar[r]^{(\ZZ /2)^3}&\PP^1 \times
\PP^1 \times \PP^1 \ar[r]^{(\ZZ/2)^3}&\PP^1
\times
\PP^1
\times
\PP^1. }
\end{equation*}
$X'$ is the normal $(\ZZ/2 \ZZ)^2$-covering of $Y$ whose resolution 
is a Burniat surface.

We have the following:

\begin{fact} \
\begin{enumerate}
  \item $\hat{X} \rightarrow X'$ is \'etale (with group $G= (\ZZ/ 2 
\ZZ)^2$) $\iff$ $S$ is a primary Burniat surface.
\item $\hat{X} \subset \sE_1 \times \sE_2 \times \sE_3$ is a 
hypersurface of multidegrees $(2,2,2)$.
\item $X' = \hat{X}/G$ is the quotient of a free action except for 
some $A_1$-singularities with stabilizer $\ZZ/ 2 \ZZ$,
yielding $\frac 14(1,1)$ - points on $X'$.
\end{enumerate}
\end{fact}

\subsection{Idea of the proof of theorem \ref{primary}}  Assume that 
$S$ is homotopically equivalent to a primary Burniat
surface. Then $S$ has the same fundamental group as a primary Burniat 
surface. Hence there is an \'etale $(\ZZ/2
\ZZ)^3$-covering $\hat{X} \rightarrow X$ with $q(\hat{X}) = 3$.
\subsection*{Step 1.} One shows that the Albanese variety 
$\Alb(\hat{X})$ is the product of three elliptic curves $\sE_1
\times \sE_2 \times
\sE_3$. In fact, for each $i \in \{1,2,3\}$ there is an intermediate cover
$$
\hat{X} \rightarrow X_i \rightarrow X'
$$ with $q(X_i) = 1$. By the universal property of the Albanese 
variety we get a morphism
$$
\lambda \colon \Alb(\hat{X}) \rightarrow \Alb(X_1) \times \Alb(X_2) 
\times \Alb(X_3) \cong \sE'_1 \times \sE'_2 \times
\sE'_3.
$$

Looking at the fundamental group $\pi_1(\hat{X})$, one sees that the 
isogeny $\lambda$ is of product type, whence the
claim follows.

\subsection*{Step 2.} Consider now the Albanese map of $\hat{X}$:
$$ f \colon \hat{X} \rightarrow f(\hat{X}) =: \hat{Y} \subset \sE_1 
\times \sE_2 \times \sE_3.
$$ Then the class of $\hat{Y}$ is the same as for the Albanese image 
of the corresponding \'etale covering of a primary
Burniat surface (since we have a map to a $ K(\pi,1)$ space with 
fundamental group $\pi$
equal to $\pi_1(\hat{X})$ hence the class of the image is invariant 
by homotopy equivalence).

Moreover, since
this class is
$2F_1 + 2F_2 + 2F_3$, we see that the Albanese map of
$\hat{X}$ is birational. Finally,  an argument using adjunction shows that
$\hat{X} \cong \hat{Y}$.

\subsection*{Step 3.}
$X$ (the canonical model of $S$) is $\hat{X}/G$, and it is a bidouble 
cover of a Del Pezzo surface of degree 6 as required.

For more details we refer to \cite{burniat1}.

\begin{rem}
  The same method proves similar theorems in the following cases:
\begin{enumerate}
  \item {\em Keum-Naie surfaces} with $K^2 = 4$ (cf. \cite{keumnaie}).
   
For Keum-Naie surfaces $S$ the key idea is to find a
representation $X =
\hat{X}/G$, where $G = (\ZZ/ 2\ZZ)^2$, acting on $\sE_1 \times 
\sE_2$ ($ \sE_1$, $ \sE_2$ being again  elliptic curves).
$\hat{X}$ is a
double cover of
$\sE_1
\times
\sE_2$ branched on a
$G$-invariant divisor $\Delta$ of bidegree $(4,4)$. 

Then $G$ acts 
freely on $\hat{X}$ for a suitable twist of the action.
\item {\em Inoue surfaces} with $K^2 = 7$ (cf. \cite{bcinoue}).

 Inoue 
surfaces are of the form $X= \hat{X}/ (\ZZ/ 2 \ZZ)^4$,
where $\hat{X}$ is a $(\ZZ/ 2 \ZZ)^4$-invariant divisor of 
multidegree $(2,2,4)$ in $\sE_1 \times \sE_2 \times D$, where
$ \sE_1$, $ \sE_2$ are again  elliptic curves, while $D$ is a curve 
of genus 5 which is a maximal abelian cover
of $\PP^1$ of exponent 2 branched on 5 points.
\item {\em Kulikov surfaces}  (cf. \cite{chancoughlan}). 

These surfaces are 
$(\ZZ / 3\ZZ)^2$-coverings of the plane branched
on the sides of a
triangle and on the three medians (the lines connecting the three 
vertices with the barycenter).
\end{enumerate}

\end{rem}

  The above cases also show another common feature:
\begin{itemize}
  \item[] $\hat{X}$ has $A_1$-singularities $\iff$ $G$ acts no longer freely.
\end{itemize} This implies that then $X = \hat{X}/G$ acquires $\nu$ 
singularities of type $\frac 14(1,1)$, and $K_S^2$ ($S$ being
a minimal model) drops by $\nu$. Furthermore $\pi_1(S) = \pi_1(X)$ 
becomes finite.

Therefore for these families  the {\em investigation of the connected 
components} of the moduli spaces has to be done by
\begin{enumerate}
  \item showing openness of the subset of the moduli space induced by 
such a family using local deformation theory;
\item investigating the closure via 1-parameter limits.
\end{enumerate}
\begin{rem}
  All the Burniat surfaces $X$ we consider are $G= (\ZZ / 2 
\ZZ)^2$-covers of a normal Del Pezzo surface $Z$ of degree
$K_X^2$.

For nodal Burniat surfaces with $K_S^2 = 4$ and Burniat surfaces with 
$K_S^2 =3$, we need to introduce a larger family,
including the
so-called {\em extended Burniat surfaces}. They will be introduced in 
the next section in the more symmetric  case of
tertiary Burniat surfaces.
\end{rem}

\section{Extended   Burniat surfaces} We recall the following 
definitions from \cite{burniat3}.
\medskip
\noindent Let $P_1, P_2, P_3 \in \PP^2$ be three non collinear 
points, and let $P_4,  \dots , P_{3+m},
\  \ m=2, 3,$ be further (distinct)  points   not lying on the sides 
of the triangle with vertices $P_1, P_2, P_3$.

Assume moreover that, for $m=2$, the points  $P_1, P_4, P_5$ are 
collinear, while, for $m=3$, we shall moreover assume
that also $P_2, P_4, P_6$ and $P_3, P_5, P_6$ are collinear (in 
particular, no four points are collinear).

Let's denote by
$\tilde{Y}:=\hat{\PP}^2(P_1, P_2,  \dots , P_{3+m})$ the weak Del 
Pezzo surface of degree $6-m$,  obtained by
blowing  up
$\PP^2$ in the $3 + m$ points
$P_1, P_2,
\dots , P_{3+m}$.

Saying that $\tilde{Y}$ is a weak Del Pezzo surface means that the 
anticanonical divisor $- K_{ \tilde{Y}}$ is nef and big;
         in our case it is not ample, because of the existence of 
(-2)-curves, i.e. curves $N_i \cong \PP^1$, with $ N_i \cdot K_{
\tilde{Y}} = 0$.

Contracting the (-2)-curves $N_i$ we obtain a normal singular Del Pezzo surface
$Y'$ with $- K_{ Y'}$ very ample.

In order to simplify the formulae, let us treat the case $m=3$, 
denoting $P'_3 : = P_4,P'_2: = P_5,P'_1: = P_6$.

Then we have that $P_i, P'_{i+1}, P'_{i+2}$ are collinear (here $i \in
\ZZ/3 \ZZ$).

We denote by $L$ the divisor on $\tilde{Y}$ which is the total 
transform of a general line in $\PP^2$, by
$E_i$ the exceptional curve lying over $P_i$,   by
$E'_i$ the exceptional curve lying over $P'_i$, and by $D_{i,1} $ the 
unique effective divisor in $ |L - E_i- E_{i+1}|$, i.e., the
proper transform of the side of the triangle joining the points $P_i, P_{i+1}$.

For $m=3$ we  have (-2)-curves $N_1, N_2, N_3$ such that
         $$\{N_i \} = | L -E_i - E'_{i+1} - E'_{i+2}|,$$

         Therefore the anticanonical image of $\tilde{Y}$ is a normal surface
$Y' \subset \PP^{6-m}$
         of degree $6-m$, whose singularities are one node $\nu_1$ (an 
$A_1$ singularity) in the case $m=2$,
         and three nodes  $\nu_1, \nu_2, \nu_3$ in the case $m=3$ (the 
(-2)-curve $N_i$ is the total transform of
         the point $\nu_i$).

         \begin{defin}\label{df}
1) The {\em Burniat branch divisors for $m=3$} are defined to be the 
divisors $D_1, D_2, D_3$ such that:
         $$
         \{D_i\}= |L - E_i- E_{i+1}| +  N_i  + |L- E_i - E'_i| +E_{i-1},
         $$
         \noindent
       2) The {\em strictly extended Burniat branch divisor classes 
for $m=3$} are defined as follows:
         $$
         \Delta_i \equiv D_i - N_i + N_{i-1} + N_{i+1},
         $$
   \noindent
       3) The {\em strictly extended Burniat branch divisors  for 
$m=3$} are defined taking an irreducible curve
$$
\Ga_i
\in |2L - E_i -E'_i-E_{i+1}-E'_{i+1}|
$$
and replacing in $\De_i$
         $$
         N_{i-1} + |L - E_i- E_{i+1}| + E_{i-1}
         $$
by $\Ga_i$, so that
         $$
         \Delta_i  = \Ga_i  + N_{i+1} + |L- E_i - E'_i|.
         $$
     \end{defin}
For each Burniat divisor, we have the option to replace it (or not) by a 
strictly extended Burniat divisor.
By taking the corresponding bidouble cover, we obtain an extended 
Burniat surface.

\begin{rem}\label{difference}

1) Observe that $(D_1 + D_2 + D_3) \in | - 3 K_{\tilde{Y}}|$ is a 
reduced normal crossing divisor.

2)  Similarly,  $( \De_1 + \De_2 + \De_3 ) \in | - 3 K_{\tilde{Y}} +
\sum N_i |$
         is a reduced normal crossing divisor.

3) On the normal Del Pezzo surface $Y'$, for $m=3$,

         $\De_j$ yields a conic and one line, $D_j$ yields three lines.

In particular, if the conic corresponding to  $\De_j$ specializes to 
contain the line corresponding to $E_{j-1}$,
         we obtain $D_2$ subtracting the divisor $ N_{j-1} + N_{j+1}$ 
and adding the divisor $N_j$.

\end{rem}

We can now consider (cf. \cite{ms}, \cite{sbc}) the associated 
bidouble covers  $S \ra \tilde{Y}$ with branch divisors the
Burniat divisors, respectively the extended Burniat divisors.

\begin{defin}
           A  {\em tertiary nodal Burniat surface}  is obtained, for $m=3$,
           as  a bidouble cover  $S \ra \tilde{Y}$ with branch 
divisors the three Burniat divisors.

  $S$ is then a minimal surface of general type with $p_g(S) =q(S)= 
0$,  $K_S^2 = 6-m$ (cf. \cite{burniat2}).

           If we let some of the three branch divisors be extended 
Burniat divisors, then we obtain a non minimal surface $S'$
whose minimal model $S$ is called
         a  {\em tertiary extended   Burniat surface}.

\end{defin}

In \cite{burniat3} it is shown that these (discontinuous) 
degenerations of the branch divisors from extended to nodal Burniat
surfaces produce a flat family of $G$-covers $X_t \rightarrow Y'$ of 
the canonical models over a normal Del Pezzo surface of
degree $6-m$ (4, respectively 3).

More precisely, we have the following two auxiliary results:
\begin{prop}\label{famiglia4} There exists a family, with connected base
$$B \subset  \{ (C_1 , \Ga_2) | C_1 \in | L - E_1| , \Ga_2 \in |2L - 
E_2 - E_3 -E_4 -E_5| \}$$ where $ C_1$ is irreducible and
either
$\Ga_2   $ is irreducible, or splits as $ N_1 + E_1 + |L-E_2 - 
E_3|$), parametrizing a flat family of canonical models, including
exactly all the nodal Burniat surfaces  and the extended  Burniat 
surfaces with $K^2_X = 4$.
\end{prop}

\begin{prop}\label{famiglia3} There exists a family, with connected base
$$T \subset  \{ (\Ga_1 , \Ga_2, \Ga_3) \}$$ where $ \Ga_1 , \Ga_2,
\Ga_3  $ are as in Definition
\ref{df}, parametrizing a flat family of canonical models, including 
exactly all the nodal Burniat surfaces  and
the extended  Burniat surfaces with $K^2_X = 3$.
\end{prop}

\begin{rem} 1) In the nodal Burniat case the surface $S$ does not 
have an ample canonical divisor $K_S$, due to the
existence of (-2)-curves, which are exactly the inverse images of the 
(-2)-curves $N_i \subset \tilde{Y}$.

For this reason we call the above Burniat surfaces
         of  {\em nodal type}. We denote their canonical model by $X$, 
and observe that
$X$ is a finite bidouble cover of the normal Del Pezzo surface $Y'$.

For $m=2$ $X$  has precisely one node (an
$A_1$-singularity, corresponding to the contraction of the 
(-2)-curve) as singularity. While, for $m=3$, $X$ has exactly three
nodes as singularities.

2) In the extended  Burniat case $S'$ is not minimal. In the strictly 
extended   Burniat case the inverse image of each $N_i$
        splits as the union of two disjoint (-1)-curves. In this 
latter case $S$ has ample canonical divisor, hence $S=X$.

3) In all cases, the morphism $ X \ra Y'$ is exactly the bicanonical 
map of $X$ (see \cite{burniat2}. \cite{burniat3}).

4) Nodal Burniat surfaces are parametrized by a family with smooth 
base of dimension $2$ for $m=2$, of dimension $1$ for
$m=3$.

         Strictly extended  Burniat surfaces are parametrized by a 
family with smooth base of dimension $3$ for $m=2$, of
dimension $4$ for
$m=3$.
\end{rem}

The key feature is that, both for nodal Burniat surfaces, and for 
extended  Burniat surfaces, the canonical model $X$ is a
finite bidouble cover of a singular Del Pezzo surface $Y'$, which has 
one node in the case $m=2$, and three nodes for $m=3$
(in this case $Y'$ is a cubic surface in $\PP^3$).

In this case the direct image $p_* (\hol_X)$ splits as a direct sum 
of four reflexive character sheaves of generic rank $1$.

  \noindent
  For $K^2 = 3$ it is shown in \cite{burniat3} that a small
deformation of a Burniat surface or of an extended Burniat surface is
  a Galois covering with group $(\mathbb{Z} / 2 \mathbb{Z})^2$ of a 
cubic surface with  three singular points,
  and with branch locus equal to three plane sections. 
  Hence one sees that  the locus of Burniat and extended Burniat surfaces is open.
   Moreover  in 
loc.cit. it is shown that the closure of the subset
corresponding to extended Burniat surfaces with $K_S^2 =3$ contains 
$G$-covers of   cubic surfaces
  with a $D_4$-singularity,
and $G$-covers of  the four nodal cubic.

Yifan Chen shows in his Bayreuth Ph.D. thesis that there are no 
further degenerations.

Summarizing, we have the following theorem
\begin{theo}[Bauer,Catanese, Chen] The irreducible component 
$\mathcal{N}$ of the moduli space containing the Burniat
surfaces with $K_S^2=3$ consists exactly of
\begin{enumerate}
  \item Burniat surfaces,
\item extended Burniat surfaces,
\item $G$-covers of a normal cubic with a $D_4$-singularity,
\item $G$-covers of the four nodal cubic, which are \'etale exactly 
over one of the four nodes.
\end{enumerate} Moreover, 1),2),3) are contained in $\mathcal{N} 
\setminus \partial \mathcal{N}$.
\end{theo}

The key technique used in the above theorem (developed in 
\cite{burniat3}) is the one of blowing up and down logarithmic
sheaves in order to calculate the tangent cohomology. It would take 
too long to explain this technique in detail here.

There remains the challenging
\begin{question}
  Is $\mathcal{N}$ a connected component of the moduli space?
\end{question}

Another approach was proposed to construct a family of surfaces 
including the tertiary Burniat surfaces, in \cite{nevespigna};
the deformation theoretic aspects were not addressed in  
\cite{nevespigna} and it  could be interesting to do it.

Using the techniques developed for the above results Yi-fan Chen has 
been able to prove a conjecture of Mendes-Lopes and
Pardini (cf. \cite{mlp2}):

\begin{theo}[Y. Chen] The six dimensional family constructed by 
Mendes -Lopes and Pardini in \cite{mlp2}, containing the
Keum-Naie surfaces with $K_S^2 =3$ as a proper algebraic subset, is 
indeed an irreducible component of the moduli space of
surfaces of general type.

\end{theo}

The new idea which made the long sought for proof of the above result 
possible is the representation of a special surface in
the above family as some
$(\ZZ / 2 \ZZ)^2$-cover of a four nodal cubic and the use of the 
methods mentioned above (to calculate
spaces of sections of
logarithmic differential forms on blow ups).
Chen shows that the tangent dimension of Kuranishi space  is bounded 
from above by 6, and then the Kuranishi
inequality  gives   that the dimension is exactly 6.

\section{Deformation of automorphisms.}

Here is what we have learnt from extended Burniat surfaces.

In this section $S$ will be the minimal model of a nodal Burniat surface with
$K_S^2 = 4$ or
$K_S^2 = 3$, and $X$  its canonical model. Observe that for $K^2 =4$, 
$X$ has one ordinary node, while for $K^2 =3$,
$X$ has three ordinary nodes.

A very surprising and new phenomenon occurs for these surfaces, 
confirming Vakil's `Murphy's law' philosophy
(\cite{murphy}).

To explain what happens for the moduli spaces of extended and nodal 
Burniat surfaces, let us  recall again an old result due
to Burns and Wahl (cf. \cite{burnswahl}).

Let $S$ be a minimal surface of general type and let $X$ be its 
canonical model. Denote by $\Def(S)$, resp. $\Def(X)$, the
base of the Kuranishi family of $S$, resp. of $X$.

Their result explains the relation between $\Def(S)$ and $\Def(X)$.

\begin{theo}[Burns - Wahl]
  Assume that $K_S$ is not ample and let $p:S \ra X$ be the  canonical morphism.

   Denote by $\mathcal{L}_X$ the space of local deformations of the 
singularities of $X$ and by
$\mathcal{L}_S$ the space of deformations of a neighbourhood of the 
exceptional curves of $p$. Then
$\Def(S)$ is realized as the fibre product associated to the Cartesian diagram

\begin{equation*}
\xymatrix{
\Def(S) \ar[d]\ar[r] & \mathcal{L}_S \cong \CC^{\nu}, \ar[d]^{\lambda} \\
\Def(X) \ar[r] & \mathcal{L}_X \cong \CC^{\nu} ,}
\end{equation*} where $\nu$ is the number of rational $(-2)$-curves in $S$, and
$\lambda$ is a Galois covering with Galois group $W := \oplus_{i=1}^r 
W_i$, the direct sum of the  Weyl groups $W_i$ of
the singular points of
$X$.
\end{theo}

An immediate consequence is the following

\begin{cor}{\bf  (Burns - Wahl)}
   1) $\psi:\Def(S) \ra \Def(X)$ is a finite morphism, in
particular, $\psi$ is surjective.

\noindent
2) If the derivative of $\Def(X) \ra \mathcal{L}_X$ is not surjective 
(i.e., the
singularities of $X$ cannot be independently
    smoothened  by the first order infinitesimal deformations of $X$),
    then $\Def(S)$ is singular.
\end{cor}

Assume now that we have $1 \neq G \leq \Aut(S) = \Aut(X)$.

Then we can consider the space of $G$-invariant local deformations of
$S$, $\Def(S,G)$, resp. $\Def(X,G)$ of  $X$, and we have a  natural 
map $\Def(S,G) \ra \Def(X,G)$.

We indeed show here that, unlike the case for  the corresponding 
morphism of   local deformation spaces, this map needs
not to be surjective; and, as far as we know, the following result 
gives the first global  example of such a  phenomenon.

\begin{theo}\label{path}
   The deformations of nodal Burniat surfaces with $K^2_S =4,3$ to 
extended  Burniat surfaces with $K^2_S =4,3$ yield
examples where
$\Def(S,(\ZZ/2\ZZ)^2) \ra \Def(X,(\ZZ/2\ZZ)^2)$  is not surjective.

Moreover, $\Def(S,(\ZZ/2\ZZ)^2) \subsetneq \Def(S)$, whereas for the 
canonical model we have: $\Def(X,(\ZZ/2\ZZ)^2) =
\Def(X)$.

The moduli space of pairs, of an  extended (or nodal) Burniat surface 
with $K^2_S =4,3$ and   a $(\ZZ/2\ZZ)^2$-action, is
disconnected; but its image in the moduli space is a connected open set.
\end{theo}

\noindent The reason for this phenomenon can already be seen locally 
around the node.

Let $G$ be the group $G \cong (\ZZ/2 )^2$ acting  on $\CC^3$ as follows:
$$G = \{1,\sigma_1, \sigma_2, \sigma_3 = \sigma_1+\sigma_2\}$$ acts by
$\sigma_1(u,v,w) = (u,v,-w)$, $\sigma_2(u,v,w) = (-u,-v,w)$.

The invariants for the action of  $G$ on $\CC^3 \times \CC$ are:
$$ x: =u^2, y:  = v^2, z : = uv , s: = w^2, t.$$

Observe that the hypersurfaces $X_t = \{ (u,v,w)| w^2 = uv + t\}$ are 
$G$-invariant,
and the quotient $X_t / G $ is the fixed hypersurface
$$ Y_t \cong Y_0 =  \{ (x,y,z)| z^2 = xy\} ,$$
which has a nodal singularity at the point $x=y=z=0$.

In fact, $G$ acts on the family $\mathcal{X}_t = \{w^2 = uv +
t\}$, which admits a simultaneous resolution only after the base 
change $\tau^2 = t$: and then we have two small
resolutions
$$
\mathcal{S}:=\{((u,v,w,\tau),\xi) \in \mathcal{X} \times \mathbb{P}^1 
| \frac{w-\tau}{u} = \frac{v}{w+\tau} = \xi\},
$$
$$
\mathcal{S'}:=\{((u,v,w,\tau),\eta) \in \mathcal{X} \times 
\mathbb{P}^1 | \frac{w+\tau}{u} = \frac{v}{w-\tau} = \eta\}.
$$ Then it is easy to see that $G$ has several liftings to $\mathcal{S}$, but
\begin{itemize}
  \item either $G$ acts only as a group of {\em birational not 
biregular} automorphisms on $\mathcal{S}$ and leaves $\tau$
fixed,
\item or, $G$ acts as a group of {\em biregular} automorphisms on 
$\mathcal{S}$ but does not leave $\tau$ fixed.
\end{itemize}

Looking at the local picture in detail, one sees how the above family 
yields   a discontinuous variation of the
three branch divisors
on the blow up $\tilde{Y_0}$ of  $Y_0$ at its singular point.
\medskip

{\bf Acknowledgement.} I would like to thank Ingrid Bauer for  helping to finish the article
without too much delay, the responsibility of misprints (I hope not too many) remaining mine.
I would like to thank prof. Katsura and  the organizers of the Conference, prof. Iritani, prof. Kawaguchi, prof. Nakaoka,  
 for their kind hospitality.


\end{document}